\documentclass[12pt]{amsart}

\usepackage{amsmath,amssymb,enumerate}
\usepackage{epsf,epsfig,amsfonts,graphicx,color}  
 \usepackage[applemac]{inputenc} 

 \date{July, 2012}   
  
 \numberwithin{equation}{section}   
\newtheorem{theorem}{\rm\bf Theorem}

\newtheorem{fact}[theorem]{\rm\bf Fact}
\theoremstyle{definition}
\newtheorem*{definition}{\rm\bf Definition}

\newtheorem{example}[theorem]{\rm\bf Example}
\newcommand{\weg}[1]{}
\newcommand{\be}{\begin{equation}}
\newcommand{\ee}{\end{equation}}

\title[Submaximal dimension of the group of almost isometries.]{On submaximal dimension of the group of almost isometries of Finsler metrics.   } \date{}
\author{Vladimir S. Matveev}\thanks{\small Institute of Mathematics,  Friedrich-Schiller-Universit\"at Jena,  07737 Jena Germany; \   vladimir.matveev@uni-jena.de \\ Partially supported by DFG (GK 1523) and FSU Jena} 

\begin{document}

\maketitle

\weg{\begin{abstract}} 
{\sc Abstract.} {\small We show that  the second greatest  possible dimension of 
the group of (local) almost isometries of a Finsler metric 
 is $\tfrac{n^2 -n}{2}  +1$ for $n= dim(M)\ne 4 $  and $\tfrac{n^2 -n}{2}  +2 =8$ for $n=4$. 
 If a Finsler metric has the  group of  almost isometries of dimension greater than $\tfrac{n^2 -n}{2}  +1$,
 then the Finsler metric  is  Randers, i.e., $F(x,y)= \sqrt{g_x(y,y)} + \tau(y)$.  Moreover,  if $n\ne 4$,    the Riemannian metric   $g$ has constant sectional  curvature and, if in addition $n\ne 2$, 
  the 1-form  $\tau$  is closed, so (locally) the metric admits the group of local isometries of the maximal 
  dimension $\tfrac{n(n+1)}2$. 
  
  In the remaining  dimensions 2 and 4, we describe all examples of Finsler  metrics with $3$ resp. 
  $8$-dimensional group of almost isometries. 
     }
\weg{\end{abstract} }

\section{Definition and results} 

Let $(M,F)$ be a connected 
 Finsler manifold of dimension $n\ge 2$. We assume that all objects in  our paper are sufficiently smooth.
   We require  that $F$ is  strictly convex but allow $F$ to be not reversible. 
 
 \begin{definition}[\cite{miguel}] 
 A  diffeomorphism $\phi:U\to V$, where $U,V\subseteq M$, 
  is called an \emph{almost isometry}, if $T(p,q,r)= T(\phi(p), \phi(q), \phi(r))$ for all $p,q,r\in U $, where $T$ is the ``triangular'' function given by 
 \begin{equation}\label{T}
 T(p,q,r):= d(p,q)+ d(q,r)- d(p,r), 
 \end{equation}
 where $d$ is the  (generally, nonsymmetric) distance corresponding to  $F$, i.e.
 $$
 d(p,q)= \textrm{inf} \left\{\int_0^1 F(c(t), c(t)') dt \mid c \in C^1([0,1]; M) \textrm{ with $c(0)= p$, $c(1)= q$}\right\}. 
 $$
 \end{definition} 
 
Similarly, a vector field $K$ is an {\it almost Killing}, if its local flow acts by almost isometries. 
 
 \begin{example} \label{Ex1}  Let $\phi:U\to V$ be  a diffeomorphism 
  such that $\phi_*(F) = F + df $, where $f  $ is a smooth function (the metric transformation $F\mapsto F+ df$ was called trivial projective change in \cite{how}). 
 Then, $\phi$ is an almost isometry.
 \end{example} 
 
 Indeed, the transformation $F \mapsto F+df$  adds $f(q) - f(p)$ (resp. $f(r)-f(q)$, resp. $f(p)- f(r)$) 
  to the length of every curve connecting $p$ and $q$ (resp. $q$ and $r$, resp. $r$ and $p$) and therefore does not change the value of $T(p,q,r)$. 
 
 \begin{fact}[\cite{miguel}, Proposition 3.2] \label{fact} 
Every almost isometry is as in Example \ref{Ex1}. \end{fact}

In our paper we answer the following natural question:
{\it what  is the submaximal dimension of the  group of almost isometries, and wenn  it is greater than the submaximal dimension of the group of isometries?} Since in view of Fact \ref{fact} every almost isometry of $F$ is an isometry of the symmetrized Finsler metric 
$F_{\textrm{sym}}(x,y)= \tfrac{1}{2}(F(x,y) + F(x,-y))$, 
  the dimension  of the group of almost isometries can not be greater than the dimension of the group of isometries of $F_{\textrm{sym}}$ which, in view of \cite[Theorem 3.3]{HD1} and \cite[\S 3]{troyanov}, is at most $\frac{n(n+1)}{2}-$dimensional. It is easy to construct the examples of (nonriemannian) reversible 
  Finsler metrics whose group of isometries (and therefore the group of almost  isometries) has dimension  $\tfrac{n^2- n}{2} +1$, see \cite[\S 3]{troyanov}.
  
    All  our  considerations are local, so we actually speak not about the group of almost isometries, but about the linear vector space of almost  Killing vector fields. Most   examples and  statements  survive  or could be generalized for  the  global setup though.

\begin{example} \label{Ex2} Let $g$ be a Riemannian metric of constant sectional 
 curvature (locally 
it has $\frac{n(n+1)}{2}-$dimensional space of Killing vector field). Consider a closed  1-form $\tau $ such that its $g$-norm  $|\tau|_g<1 $ at all points. Then, any Killing vector field of $g$ is an  almost Killing vector field of the Randers Finsler metrics $F_{g, \tau}$ defined  by $F_{g,\tau}(x,y)= \sqrt{g_x(y,y)} + \tau(y) , $ and vice versa.
\end{example} 
Indeed, the implication ``$\Rightarrow$'' follows from Example \ref{Ex1}  and the implication ``$\Leftarrow$'' follows from the above mentioned observation that almost isometries are isometries of the symmetrized metric (which is the case of Randers metrics  is  essentially the initial metric $g$).  

The main results of our papers are  Theorems \ref{thm1}, \ref{thm2}, \ref{thm2.5}, \ref{thm3} and Examples \ref{Ex2.5}, \ref{Ex3} below. 

\begin{theorem} \label{thm1}
Assume the vector space of almost Killing vector fields on $(M,F)$
is  more than $ \tfrac{n^2- n}{2}  +1$- dimensional. Then, $F$ is a Randers metric, i.e., 
$F(x,y)=   \sqrt{g_x(y,y)} + \tau(y)$ for a  Riemannian metric $g$ and for a 1-form $\tau$ with $|\tau|_g<1$. 
\end{theorem}

\begin{theorem} \label{thm2} 
  Let $F= F_{g,\tau}$ be a  Randers Finsler metric. Assume $n= dim(M)\ne 2,4$.  Suppose the  space of almost Killing vector fields is  more than $ \tfrac{n^2- n}{2}  +1$- dimensional. 
  Then, $F$ is as in Example \ref{Ex2}, i.e., $g$ has constant sectional curvature and the form $\tau $ is closed.  
\end{theorem}

The above two theorems answer our question for all dimensions except of $n=2,4$, and, in a certain sense, tell us  that  for metrics with many Killing resp. almost Killing vector fields there is no big difference between isometries and almost isometries.  The case 
 $n=2$ will be considered in Example \ref{Ex2.5} and  Theorem \ref{thm2.5}; we will see that in this case the submaximal dimension of 
 the space of almost Killing vector fields  is still $\tfrac{n^2 - n}{2} +1= 2$, but the description of  Finsler metrics with the space of almost Killing vector fields  of the  highest  dimension $3$   is slightly more complicated (the metric is still of constant curvature but the form $\tau $  may be  not closed).  
  
  The  4-dimensional case  (considered in  Example \ref{Ex3} and 
  Theorem \ref{thm3}) is much  more interesting: 
  remarkably, there exist   examples  of  Randers Finsler metrics such that they are 
  not as in Example \ref{Ex2} and such that the dimension of the space of almost Killing vector fields is $8$ which is greater  than $\tfrac{n^2- n}{2}  +1= 7$; we construct them all. 

\begin{example} \label{Ex2.5} Consider a Riemannian $(n=2)$-dimensional manifold $(M^2, g)$  
of constant sectional curvature. Let $\omega= Vol_g = \sqrt{det(g)} dx\wedge dy$ be the volume form of $g$ and $\tau$ be a 1-form such that for a certain constant $c$ we have $d\tau= c \cdot \omega$. Assume that $|\tau|_g <1$    so   $F_{g,\tau}$ is a Finsler metric. Then, locally,  the space of almost Killing vector fields for this metric is   $\tfrac{n^2- n}{2}  +1+1= 3$-dimensional. 
\end{example}

Note that  Example \ref{Ex2.5}  is essentially local and can not live 
 on closed manifolds, since the volume  form  on a closed manifold 
 can not be differential of a 1-form by  Stokes' theorem.  
\begin{example} \label{Ex3}   Consider a K\"a{}hler 4-manifold $(M^4,g, \omega)$ of constant holomophic curvature. Take a  1-form $\tau$ such that  $d\tau= c \cdot \omega$ for a certain  constant $c$.    Assume that $|\tau|_g <1$    so   $F_{g,\tau}$ is a Finsler metric.  Then, locally, 
the space of almost Killing vector fields for this metric is  $\tfrac{n^2- n}{2}  +1+1= 8$-dimensional. 
\end{example}

Indeed, 
consider the space of Killing vector fields preserving the volume form $\omega$ in dimension $2$ or 
 the K\"ahler symplectic form $\omega$ in dimension $4$.  This space is 3-dimensional in dimension 2 (since every orientation-preserving isometry preserves the volume form)  and 
   $8$-dimensional in dimension 4 (since  a K\"ahler space of constant holomorphic 
    curvature is  a locally symmetric space with $\omega$-preserving-isotropy  subgroup isomorphic  to the  4-dimensional $U_2$). 
     Every such Killing vector field is an almost Killing vector field of $F_{g,\tau}$. 
Indeed, its  flow is an isometry of $g$  and 
 sends $\tau $ to another 1-form such that its differential is still $c \cdot \omega$ (because this  isometry  preserves $\omega$).  

\begin{theorem} \label{thm2.5}  Suppose a Randers Finsler metric is not as in Example \ref{Ex2}. Assume $n= dim(M)= 2$. 
Then, if the  space of almost Killing vector fields for this metric is (at least)  $\tfrac{n^2- n}{2}  +1+1= 3$-dimensional, then the metric is as in   Example \ref{Ex2.5} (with $c\ne 0$).
\end{theorem} 

\begin{theorem} \label{thm3} Suppose a Randers Finsler metric is not as in Example \ref{Ex2}. Assume $n= dim(M)= 4$. 
Then, if the  space of almost Killing vector fields for this metric is (at least)  $\tfrac{n^2- n}{2}  +1+1= 8$-dimensional, then the metric is as in Example \ref{Ex3}. If the space of almost Killing vector fields has a higher dimension, the metric is as in Example \ref{Ex2}.
\end{theorem}

\section{Proofs}

{\bf Proof of Theorem \ref{thm1}.} 
 Assume a metric $F$ on $M^n$ has at least $\frac{n^2-n}{2}+2$ linearly independent
 almost Killing vector fields. Let us slightly improve the metric $F$, i.e., construct canonically a ``better'' metric $F_{better}$ such that each almost Killing vector field of $F$ is a Killing vector field of $F_{better}$. This ``improvement'' will take place in each tangent space independently, though of course the result will still smoothly depend on the point.  Take a point $p\in M $
 and consider $T_pM$. Consider the unit ball $K_F$ in the norm $F_{|T_pM}$:
 $$
 K_F:= \{y \in T_pM \mid F(p,y)\le 1\}.
 $$ 
 It is a convex body in $T_pM$ containing the zero vector $\vec 0\in T_pM$.
 Now, consider the dual space $ T_p^*M$ and  the dual(=polar) convex body there:
 $$
 K_F^*:= \{ \xi\in T^*_pM \mid \xi(y)\le 1 \textrm{ for all $y \in K_F$} \}.$$
 Take the barycenter $b_F$ of $K_F^*$ and consider   the convex body $$ K^*_{better}:= K_F^*- b_F = \{\xi - b_F\mid \xi \in K_F^*\}.$$ The barycenter of $K^*_{better}$ lies at $\vec 0$. 
  
  Next, consider the body $K_{better}\subseteq T_pM$  dual to $K^*_{better}$ and the Finsler metric $F_{better}$ such that at every point $p$ the unit ball is the corresponding $K_{better}$. 
  Evidently, $F_{better}$ is a smooth Finsler metric. 

\begin{example} \label{Ex4} Let $F(x,y)= \sqrt{g_x(y,y)} + \tau(y)$  be a Randers metric. Then, $F_{better}(x,y) = \sqrt{g_x(y,y)} $, i.e., is essentially the  Riemannian metric $g$.  Moreover, if $F_{better}$ is essentially a Riemannian metric, then $F$ is a Randers metric (to see all this, it is recommended to consider a $g$-orthonormal basis in $T_xM$ and calculate everything in the corresponding coordinate system, which is an easy exercise). 
   \end{example} 

Let us now show that the transformation $F \mapsto F+ \tau$ (where $\tau $ is a 1-form) does not change the metric $F_{better}  $, in the sense that the metrics $F_{better}$ constructed for $F$ and for $F + \tau$ coincide. This fact is well known in the convex geometry, we prove it for convenience of the reader.  

First, since the sets $K_F$, $K_{F+\tau}$ are convex and compact, the maximum of every 1-form $\xi\in T^*_pM$ over $K_F$ resp.  $K_{F+\tau}$ is achieved at a point of the unit sphere  $S_F= \{y \in T_pM\mid F(p,y)= 1\}$ resp. $S_{F+\tau}= \{y \in T_pM\mid F(p,y) + \tau(y)= 1\}$. 

Consider the bijection 
$$
f:S_F\to S_{F+\tau}, \ \ f(y)= \tfrac{1}{1+ \tau(y)} y  .$$ 

Consider $\xi\in K^*_F$. Then, for every element $ y'= f(y)$ of 
$S_{F+ \tau}$   (for $y\in S_F$) we  have 
$$
(\xi + \tau)(y')=( \xi + \tau )\left(\tfrac{1}{1+ \tau(y)} y\right)= \tfrac{1}{1+ \tau(y)}\xi(y)+ \tfrac{1}{1+ \tau(y)}\tau(y)\le \tfrac{1 +\tau(y)}{1+ \tau(y)}= 1.$$

We see that for every 1-form   $\xi \in K^*_{F}$ the 1-form $\xi + \tau\in K^*_{F+\tau}$. Analogous we prove that  for every $\xi \in K^*_{F+\tau}$ the 1-form  $\xi -  \tau$ lies in $K^*_{F}$. Thus,   $K^*_{F+ \tau}$ is the  $\tau$-parallel translation of $K^*_F$. Then, the 
   barycenter  $b_{F+\tau}$ 
   corresponding to $K^*_{F+ \tau} $ 
  is $b_F + \tau$, and the bodies  
  $K_{better}$ corresponding to $F+ \tau$ and  $F$   coincide, so the transformation $F\mapsto F+ \tau$ does not change the metric $F_{better}$.

  Since by Fact \ref{fact}
   every almost isometry sends $F$  to a $F+ df$, and since the addition of    $\tau = df$ does not change the metric $F_{better}$, every almost isometry of $F$ is an isometry of $F_{better}$.

   By assumptions, the initial Finsler metric $F$ has at least $\tfrac{n^2 -n}{2} +2$-dimensional space of almost Killing vector field.
   Thus,  the Finsler metric $F_{better}$ has at least $\tfrac{n^2 -n}{2} +2$-dimensional space of Killing vector field.    By \cite[Theorem 3.1]{troyanov},  it is essentially a Riemannian metric, i.e., $F(x,y)= \sqrt{g_{x}(y,y)}$ for a certain Riemannian metric $g$. Then, 
   as we explainen in Example \ref{Ex4}, the metric $F$ is a Randers metric as we claimed.    Theorem \ref{thm1} is proved.

{\bf Proof of Theorem \ref{thm2}}. We assume $n =dim(M)\ne 2,4$ and   consider a Finsler metric $F$ 
 such that 
the vector space of its almost Killing vector fields is at least $\frac{n^2- n}{2} + 2$-dimensional. By Theorem \ref{thm1}, the  Finsler metric is a Randers one,  $F(x,y)= \sqrt{g_x(y,y)} + \tau(y)$.

Take a point $x\in M$ and consider the   Lie subgroup of $SO(T_xM; g_x)$ (=the group of $g_x$-isometries of $T_xM$ preserving the orientation)  
 generated by almost Killing vector fields that vanish at $0$. 
 It is at least  $\left(\frac{n^2- n}{2} + 2- n\right)$-dimensional.  Take it closure; it is a closed Lie subgroup of  $SO(T_xM; g_x)$ of dimension at least $\frac{n^2- n}{2} + 2- n = \tfrac{(n-1)(n-2)}{2} +1$.

By the classical result of \cite{montgomery}, for $n\ne 4$, every closed subgroup of $SO_n$ of dimension greater than $\tfrac{(n-1)(n-2)}{2}$ coincides with the whole $SO_n$. Evidently, the flow of every almost Killing vector field preserves $d\tau$ (since it preserves the metric $g$ and  in view of Fact \ref{fact} sends  
$\tau$ to  $\tau + df$ which  have the same differential as $\tau$). Then,  the differential of $\tau$, which is a 2-form,  is preserved by the whole group $SO(T_xM; g_x)$ implying it is vanishes  at the  (arbitrary) point $x$. 

Since $SO(T_xM; g_x)$  acts transitively on 2-planes in $T_xM$,   the Riemannian metric 
$g$ has constant sectional curvature. 
   Theorem 2 is proved.

  {\bf Proof of Theorem \ref{thm2.5}.}    Consider a Randers metric $F(x,y)=\sqrt{g_x(y,y)} + \tau$;  we assume $n= dim(M)= 2$  and the existence of $3$ linearly independent 
almost Killing vector fields. Then, the Riemannian metric $g$ has 3 Killing vector fields implying it has constant sectional curvature. The Killing vector fields preserve the differential of $\tau$ implying $d\tau$ is proportional to the volume form. Theorem \ref{thm2.5} is proved.   
  
{\bf   Proof of Theorem \ref{thm3}.}  Consider a Randers metric $F(x,y)=\sqrt{g_x(y,y)} + \tau$; 
 we assume $n= dim(M)= 4$  and the existence of at least $8$ linearly independent 
almost Killing vector fields.   Since almost Killing vector fields for $F$ are Killing for $g$, the  metric $g$ has at least  $8$ linearly independent Killing vector fields. 4-dimensional 
Riemannian metrics with at least 8 linearly independent Killing vector fields are all known
(see for example  \cite{ishihara}): they are 
of constant sectional curvature, or of constant holomophic sectional curvature.      In both   cases the local pseudogroup of almost isometries generated by the almost Killing vector fields acts transitively.

Suppose now $d\tau\ne 0$ at a certain point. Since   the local pseudogroup of almost isometries generated by the almost Killing vector fields acts transitively, $d\tau\ne 0$  at all points. 
 Consider an arbitrary  point $x$ 
 and  again, as in the proof of Theorem \ref{thm2},   consider the  
  almost Killing vector fields that vanish at $x$, the Lie subgroup of $SO(T_xM; g_x)$ generated by these  vector fields, and its closure which we denote $H_x\subseteq SO(T_xM; g_x) \approx SO_4$. 
   It is at least 4 dimensional. Now, the closed subgroups of $SO(T_xM; g_x)   \approx SO_4$ of dimensions $4$ and larger are well-understood: any 4-dimensional closed (connected)  subgroup  is essentially $U_2$ 
   (in the sense that in a certain $g_x$-orthonormal  basis of $T_xM$   the group 
 $  SO(T_xM; g_x) $ and its subgroup $H$ are precisely the  standard  $SO_4$ and the standard 
  $U_2$ standardly embedded in $SO_4$).  Any closed 
  (connected) subgroup  of dimension $\ge 5$  is the whole $SO_4$.

In all cases  the element $-id:T_xM\to T_xM$  is an element of $H_x$ which implies that 
every geodesic reflection 
(i.e., a local  isometry $I$ of $g$ such that it takes the point $x$ to itself and whose  differential at $x$ is $-id$) 
can be realized by an almost isometry of $F$ generated by almost 
Killing vector fields.

By our assumptions,   the differential $d\tau$ 
of the form $\tau$ is not zero. Since it is preserved by $H_x$, it also preserved by the holonomy group of $g$, since 
the  holonomy group of a symmetric space  is the subgroup of the group generated by  geodesic reflections. Then, $d\tau$ is covariantly constant. Since the group $H_x$ acts transitively on $T_xM$,   
 all  the eigenvalues 
of the endomorphism $J:TM\to TM$ 
defined by the condition  $d\tau(y, .) = g(J(y), .)  $  are $\pm c  \cdot  i$, so the endomorphism 
$\tfrac{1}{c}\cdot J$ is an almost 
 complex structure.  It is covariantly constant (because $g$ and $d\tau$ are covariantly constant), so $c$ is a constant.  Then, it is a complex structure  and $(M, g,  \omega= \tfrac{1}{c} d\tau)$ is a K\"ahler manifold.   Then, the metric has constant holomorphic sectional curvature  and  $d\tau$ is proportional to $\omega$ with the constant coefficient as we claimed.   Finally, if there exists more than 8 almost Killing vector fields, then the group $H_x$ is more than 4-dimensional, which, as we explained above, implies that it is the whole $SO(T_xM;g_x)$ so the metric has constant sectional curvature and the form $\tau$ is closed.

\paragraph{\bf Acknowledgement.}
This result was obtained and the paper was essentially written during XVII Escola de Geometria Diferencial at Manaus. My  interest to the problem is motivated by the talk of Miguel Angel Javaloyes and the paper answers  the problems  he stated  during his talk and in the discussion afterward.   I thank the organizers of the conference  for partial financial support and for their hospitality,  and 
the Deutsche Forschungsgemeinschaft (GK 1523) and FSU Jena  for  partial financial support.


\begin{thebibliography}{99}

 
\bibitem{HD1} 
S. Deng and  Z. Hou,  \emph{
The group of isometries of a Finsler space, }
Pacific J. Math. {\bf 207}(2002), no. 1, 149--155. 
 
 \bibitem{miguel}
 M. A. Javaloyes, L. Lichtenfelz, P. Piccione, \emph{ Almost isometries of non-reversible metrics with applications to stationary spacetimes,} 	arXiv:1205.4539v1 [math.DG].

 \bibitem{ishihara} S. Ishihara, 
\emph{Homogeneous Riemannian spaces,}
J. Math. Soc. Japan {\bf 7}(1955), 345--370.
 
\bibitem{Kobayashi72} S. Kobayash and T.  Nagano,
\emph{Riemannian manifolds with abundant isometries,} Differential geometry (in honor of Kentaro Yano). pp. 195--219. Kinokuniya, Tokyo, 1972.

 
 
 
\bibitem{how} V. S. Matveev,
\emph{Can we make a Finsler metric complete by a trivial projective change?} Accepted to proceedings of the VI International Meeting on Lorentzian Geometry (Granada, September 6--9, 2011), 	arXiv:1112.5060.
\bibitem{sigma}  V. S. Matveev,   	\emph{ On the dimension of the group of projective transformations,} 
 SIGMA, 8(2012), 007, 4 pages.




\bibitem{troyanov} V. S. Matveev, M. Troyanov,  \emph{The Binet-Legendre Metric  in Finsler Geometry. }  Geom \& Topol., to appear; arXiv:1104.1647. 


 \bibitem{montgomery} D. Montgomery and H. Samelson, 
  \emph{Transformation groups of spheres}, Ann. of Math. (2) {\bf 44}(1943),  454--470.
  
\weg{\bibitem{Yano} K. Yano,  \emph{On $n$-dimensional Riemannian spaces admitting a group of motions of order $\frac{n(n-1)}{2}+1$,} Trans. Amer. Math. Soc. {\bf 74}(1953), 260--279.}
\end{thebibliography}
\end{document}